\documentclass[11pt]{article}
\usepackage{graphicx}
\usepackage[round]{natbib}
\usepackage[colorlinks=true,citecolor=black,linkcolor=black,urlcolor=black]{hyperref}
\newcommand{\Prob}[1]{\ensuremath{\mathrm{Pr}\left( #1 \right)}}
\newcommand{\Expect}[1]{\ensuremath{\mathbf{E}\left[ #1 \right]}}
\begin{document}
\title{Philosophy and the practice of Bayesian statistics}
\author{Andrew Gelman\\
{\small Department of Statistics and Department of Political Science, Columbia University}\\
\and
Cosma Rohilla Shalizi\\
{\small Statistics Department, Carnegie Mellon University}\\
{\small Santa Fe Institute}}
\date{27 June 2011}
\maketitle

\begin{abstract}
  A substantial school in the philosophy of science identifies Bayesian
  inference with inductive inference and even rationality as such, and seems to
  be strengthened by the rise and practical success of Bayesian statistics.  We
  argue that the most successful forms of Bayesian statistics do not actually
  support that particular philosophy but rather accord much better with
  sophisticated forms of hypothetico-deductivism.  We examine the actual role
  played by prior distributions in Bayesian models, and the crucial aspects of
  model checking and model revision, which fall outside the scope of Bayesian
  confirmation theory.  We draw on the literature on the consistency of
  Bayesian updating and also on our experience of applied work in social
  science.

  Clarity about these matters should benefit not just philosophy of science,
  but also statistical practice.  At best, the inductivist view has encouraged
  researchers to fit and compare models without checking them; at worst,
  theorists have actively discouraged practitioners from performing model
  checking because it does not fit into their framework.
\end{abstract}

\section{The usual story---which we don't like}

\begin{quotation}
\noindent
  {\em In so far as I have a coherent philosophy of statistics, I hope it is
  ``robust'' enough to cope in principle with the whole of statistics, and
  sufficiently undogmatic not to imply that all those who may think rather
  differently from me are necessarily stupid.  If at times I do seem dogmatic,
  it is because it is convenient to give my own views as unequivocally as
  possible.}  \citep[p.\ 458]{Bartlett-inference}
\end{quotation}

Schools of statistical inference are sometimes linked to approaches to the
philosophy of science.  ``Classical'' statistics---as exemplified by Fisher's
$p$-values, Neyman-Pearson hypothesis tests, and Neyman's confidence
intervals---is associated with the hypothetico-deductive and falsificationist
view of science.  Scientists devise hypotheses, deduce implications for
observations from them, and test those implications.  Scientific hypotheses can
be rejected (that is, falsified), but never really established or accepted in
the same way.  \citet{Mayo-error} presents the leading contemporary statement
of this view.

\begin{figure}
  \begin{center}
    \resizebox{\textwidth}{!}{\includegraphics{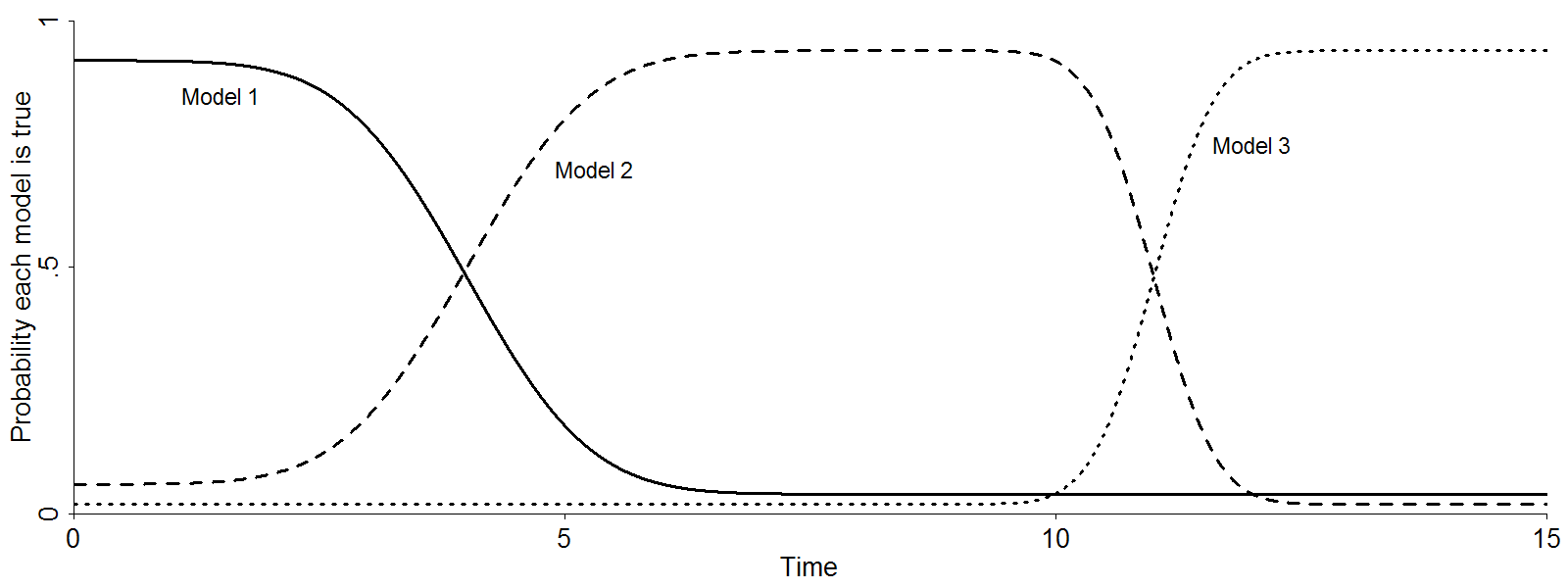}}
\vspace{-.4in}
  \end{center}
  \caption{\em Hypothetical picture of idealized Bayesian inference under the conventional inductive philosophy. The posterior probability of different models changes over time with the expansion of the likelihood as more data are entered into the analysis.  Depending on the context of the problem, the time scale on the $x$-axis might be hours, years, or decades, in any case long enough for information to be gathered and analyzed that first knocks out hypothesis 1 in favor of hypothesis 2, which in turn is
dethroned in favor of the current champion, model 3.}
  \label{fig:idealizedcurve}
\end{figure}

In contrast, Bayesian statistics or ``inverse probability''---starting with a
prior distribution, getting data, and moving to the posterior distribution---is
associated with an inductive approach of learning about the general from
particulars.  Rather than testing and attempted falsification, learning
proceeds more smoothly: an accretion of evidence is summarized by a posterior
distribution, and scientific process is associated with the rise and fall in
the posterior probabilities of various models; see Figure \ref{fig:idealizedcurve} for a schematic illustration.  In this view, the expression
$p(\theta|y)$ says it all, and the central goal of Bayesian inference is
computing the posterior probabilities of hypotheses.  Anything not contained in
the posterior distribution $p(\theta|y)$ is simply irrelevant, and it would be
irrational (or incoherent) to attempt falsification, unless that somehow shows
up in the posterior.  The goal is to learn about general laws, as expressed in
the probability that one model or another is correct.  This view, strongly
influenced by \citet{Savage-foundations}, is widespread and influential in the
philosophy of science (especially in the form of Bayesian confirmation theory;
see \citealt{Howson-Urbach,Earman-on-Bayes}) and among Bayesian statisticians
\citep{Bernardo-et-al-bayesian-theory}.  Many people see support for this view
in the rising use of Bayesian methods in applied statistical work over the last
few decades.\footnote{Consider the current (9 June 2010) state of the Wikipedia
  article on Bayesian inference, which {\em begins} as follows:
  \begin{quotation}\noindent
    Bayesian inference is statistical inference in which evidence or
    observations are used to update or to newly infer the probability that a
    hypothesis may be true.
  \end{quotation}
  It then continues with:
  \begin{quotation}\noindent
    Bayesian inference uses aspects of the scientific method, which involves
    collecting evidence that is meant to be consistent or inconsistent with a
    given hypothesis. As evidence accumulates, the degree of belief in a
    hypothesis ought to change. With enough evidence, it should become very
    high or very low.\,\dots
    Bayesian inference uses a numerical estimate of the degree of belief in a
    hypothesis before evidence has been observed and calculates a numerical
    estimate of the degree of belief in the hypothesis after evidence has been
    observed.\,\dots Bayesian inference usually relies on degrees of belief, or
    subjective probabilities, in the induction process and does not necessarily
    claim to provide an objective method of induction. Nonetheless, some
    Bayesian statisticians believe probabilities can have an objective value
    and therefore Bayesian inference can provide an objective method of
    induction.
  \end{quotation}
  These views differ from those of, e.g.,
  \citet{Bernardo-et-al-bayesian-theory} or \citet{Howson-Urbach} only in the omission of technical details. \label{note:wikipedia}}

We think most of this received view of Bayesian inference is wrong.\footnote{We are claiming that most of the standard philosophy of Bayes is wrong, {\em not} that most of Bayesian inference itself is wrong.  A statistical method can be useful even if its philosophical justification is in error.  It is precisely because we believe in the importance and utility of Bayesian inference that we are interested in clarifying its foundations.} Bayesian
methods are no more inductive than any other mode of statistical inference.  Bayesian data analysis is much
better understood from a hypothetico-deductive perspective.\footnote{We are not
  interested in the hypothetico-deductive ``confirmation theory'' prominent in
  philosophy of science from the 1950s through the 1970s, and linked to the
  name of \citet{Hempel-explanation}.  The hypothetico-deductive account of
  scientific method to which we appeal is distinct from, and much older than,
  this particular sub-branch of confirmation theory.}  Implicit in the best
Bayesian practice is a stance that has much in common with the
error-statistical approach of \citet{Mayo-error}, despite the latter's
frequentist orientation.  Indeed, crucial parts of Bayesian data analysis, such
as model checking, can be understood as ``error probes'' in Mayo's sense.

We proceed by a combination of examining concrete cases of Bayesian data
analysis in empirical social science research, and theoretical results on the
consistency and convergence of Bayesian updating.  Social-scientific data
analysis is especially salient for our purposes because there is general
agreement that, in this domain, all models in use are wrong---not merely
falsifiable, but actually false.  With enough data---and often only a fairly
moderate amount---any analyst could reject any model now in use to any desired
level of confidence.  Model fitting is nonetheless a valuable activity, and
indeed the crux of data analysis.  To understand why this is so, we need to
examine how models are built, fitted, used, and checked, and the effects of
misspecification on models.

Our perspective is not new; in methods and also in philosophy we follow
statisticians such as
\citet{Box-sampling-and-Bayes,Box-apology-for-ecumenism,Box-on-Shafer},
\citet{Good-Crook-compromise}, \citet{Good-on-Good},
\citet{Morris-on-Efron-on-Bayes}, \citet{Hill-model-based-statistics}, and
\citet{Jaynes-book}.  All these writers emphasized the value of model checking
and frequency evaluation as guidelines for Bayesian inference (or, to look at
it another way, the value of Bayesian inference as an approach for obtaining
statistical methods with good frequency peroperties; see
\citealt{Rubin-bayesianly-justifiable}).  Despite this literature, and despite
the strong thread of model checking in applied statistics, this philosophy of
Box and others remains a minority view that is much less popular than the idea
of Bayes being used to update the probabilities of different candidate models
being true (as can be seen, for example, by the Wikipedia snippets given in
footnote \ref{note:wikipedia}).

A puzzle then arises: The evidently successful methods of modeling and model
checking (associated with Box, Rubin, and others) seems out of step with the
accepted view of Bayesian inference as inductive reasoning (what we call here
``the usual story'')?  How can we understand this disjunction?  One possibility
(perhaps held by the authors of the Wikipedia article) is that the
inductive-Bayes philosophy is correct and that the model-building approach of
Box et al.\ can, with care, be interpreted in that way.  Another possibility is
that the approach characterized by Bayesian model checking and continuous model
expansion could be improved by moving to a fully-Bayesian approach centering on
the posterior probabilities of competing models.  A third possibility, which we
advocate, is that Box, Rubin, et al.\ are correct and that the usual
philosophical story of Bayes as inductive inference is faulty.

We are interested in philosophy and think it is important for statistical
practice---if nothing else, we believe that strictures derived from philosophy
can inhibit research progress.\footnote{For example, we have more than once
  encountered Bayesian statisticians who had no interest in assessing the fit
  of their models to data because they felt that Bayesian models were by
  definition subjective, and thus neither could nor should be tested.}  That
said, we are statisticians, not philosophers, and we recognize that our
coverage of the philosophical literature will be incomplete.  In this
presentation, we focus on the classical ideas of Popper and Kuhn, partly
because of their influence in the general scientific culture and partly because
they represent certain attitudes which we believe are important in
understanding the dynamic process of statistical modeling.  We also emphasize
the work of \citet{Mayo-error} and \citet{Mayo-Spanos-post-data} because of its
relevance to our discussion of model checking.  We hope and anticipate that
others can expand the links to other modern strands of philosophy of science
such as
\citet{Giere-explaining,Haack-evidence-and-inquiry,Kitcher-advancement,Laudan-beyond}
which are relevant to the freewheeling world of practical statistics; our goal
here is to demonstrate a possible Bayesian philosophy that goes beyond the
usual inductivism and can better match Bayesian practice as we know it.

\section{The data-analysis cycle}
\label{sec:data-analysis-cycle}

We begin with a very brief reminder of how statistical models are built and
used in data analysis, following \citet{Gelman-et-al-Bayesian-data-analysis},
or, from a frequentist perspective, \citet{Guttorp-scientific}.

The statistician begins with a model that stochastically generates all the data
$y$, whose joint distribution is specified as a function of a vector of
parameters $\theta$ from a space $\Theta$ (which may, in the case of some
so-called non-parametric models, be infinite dimensional).  This joint
distribution is the likelihood function.  The stochastic model may involve
other, unmeasured but potentially observable variables $\tilde{y}$---that is,
missing or latent data---and more-or-less fixed aspects of the data-generating
process as covariates.  For both Bayesians and frequentists, the joint
distribution of $(y,\tilde{y})$ depends on $\theta$.  Bayesians insist on a
full joint distribution, embracing observables, latent variables, and
parameters, so that the likelihood function becomes a conditional probability
density, $p(y|\theta)$.  In designing the stochastic process for
$(y,\tilde{y})$, the goal is to represent the systematic relationships between
the variables and between the variables and the parameters, and as well as to
represent the noisy (contingent, accidental, irreproducible) aspects of the
data stochastically.  Against the desire for accurate representation one must
balance conceptual, mathematical and computational tractability.  Some
parameters thus have fairly concrete real-world referents, such as the famous
(in statistics) survey of the rat population of Baltimore
\citep{Brown-Sallow-Davis-Cochran}.  Others, however, will reflect the
specification as a mathematical object more than the reality being
modeled---$t$-distributions are sometimes used to model heavy-tailed
observational noise, with the number of degrees of freedom for the $t$
representing the shape of the distribution; few statisticians would take this
as realistically as the number of rats.

Bayesian modeling, as mentioned, requires a joint distribution for
$(y,\tilde{y},\theta)$, which is conveniently factored (without loss of
generality) into a prior distribution for the parameters, $p(\theta)$, and the
complete-data likelihood, $p(y,\tilde{y}|\theta)$, so that $p(y|\theta) =
\int{p(y,\tilde{y}|\theta) d\tilde{y}}$.  The prior distribution is, as we will
see, really part of the model.  In practice, the various parts of the model
have functional forms picked by a mix of substantive knowledge, scientific
conjectures, statistical properties, analytical convenience, disciplinary
tradition, and computational tractability.

Having completed the specification, the Bayesian analyst calculates the
posterior distribution $p(\theta|y)$; it is so that this quantity makes sense
that the observed $y$ and the parameters $\theta$ must have a joint
distribution.  The rise of Bayesian methods in applications has rested on
finding new ways of to actually carry through this calculation, even if only
approximately, notably by adopting Markov chain Monte Carlo methods, originally
developed in statistical physics to evaluate high-dimensional integrals
\citep{Metropolis-et-al-Monte-Carlo,MEJN-on-Monte-Carlo}, to sample from the
posterior distribution.  The natural counterpart of this stage for non-Bayesian
analyses are various forms of point and interval estimation to identify the set
of values of $\theta$ that are consistent with the data $y$.

According to the view we sketched above, data analysis basically ends with the
calculation of the posterior $p(\theta|y)$.  At most, this might be elaborated
by partitioning $\Theta$ into a set of models or hypotheses, $\Theta_1, ...\Theta_K$, each with a prior probability
$p(\Theta_k)$ and its own set of parameters $\theta_k$.  One would then compute
the posterior parameter distribution within each model, $p(\theta_k|y,\Theta_k)$,
and the posterior probabilities of the models,
\begin{eqnarray*}
p(\Theta_k|y)  & = & \frac{p(\Theta_k)p(y|\Theta_k)}{\sum_{{k}^{\prime}}{(p(\Theta_{{k}^{\prime}})p(y|\Theta_{{k}^{\prime}}))}}\\
& = & \frac{p(\Theta_k)\int{p(y,\theta_k|\Theta_k)d\theta_k}}{\sum_{{k}^{\prime}}{ (p(\Theta_{{k}^{\prime}})\int p(y,\theta_k|\Theta_{{k}^{\prime}})d\theta_{{k}^{\prime}})}}.
\end{eqnarray*}
These posterior probabilities of hypotheses can be used for Bayesian model
selection or Bayesian model averaging (topics to which we return below).
Scientific progress, in this view, consists of gathering data---perhaps
through well-designed experiments, designed to distinguish among interesting
competing scientific hypotheses (cf.\
\citealp{Atkinson-Donev,Paninski-information-theoretic-experimental-design})---and then plotting the $p(\Theta_k|y)$'s over time and watching the system learn (as sketched in Figure \ref{fig:idealizedcurve}).

In our view, the account of the last paragraph is crucially mistaken.  The
data-analysis process---Bayesian or otherwise---does not end with calculating
parameter estimates or posterior distribution.  Rather, the model can then be
{\em checked}, by comparing the implications of the fitted model to the
empirical evidence.  One asks questions like, Do simulations from the fitted
model resemble the original data?  Is the fitted model consistent with other
data not used in the fitting of the model?  Do variables that the model says
are noise (``error terms'') in fact display readily-detectable patterns?
Discrepancies between the model and data can be used to learn about the ways in
which the model is inadequate for the scientific purposes at hand, and thus to
motivate expansions and changes to the model (\S \ref{sec:model-checking}).

\subsection{Example: Estimating voting patterns in subsets of the population}

\begin{figure}
  \begin{center}
    \includegraphics[width=\textwidth]{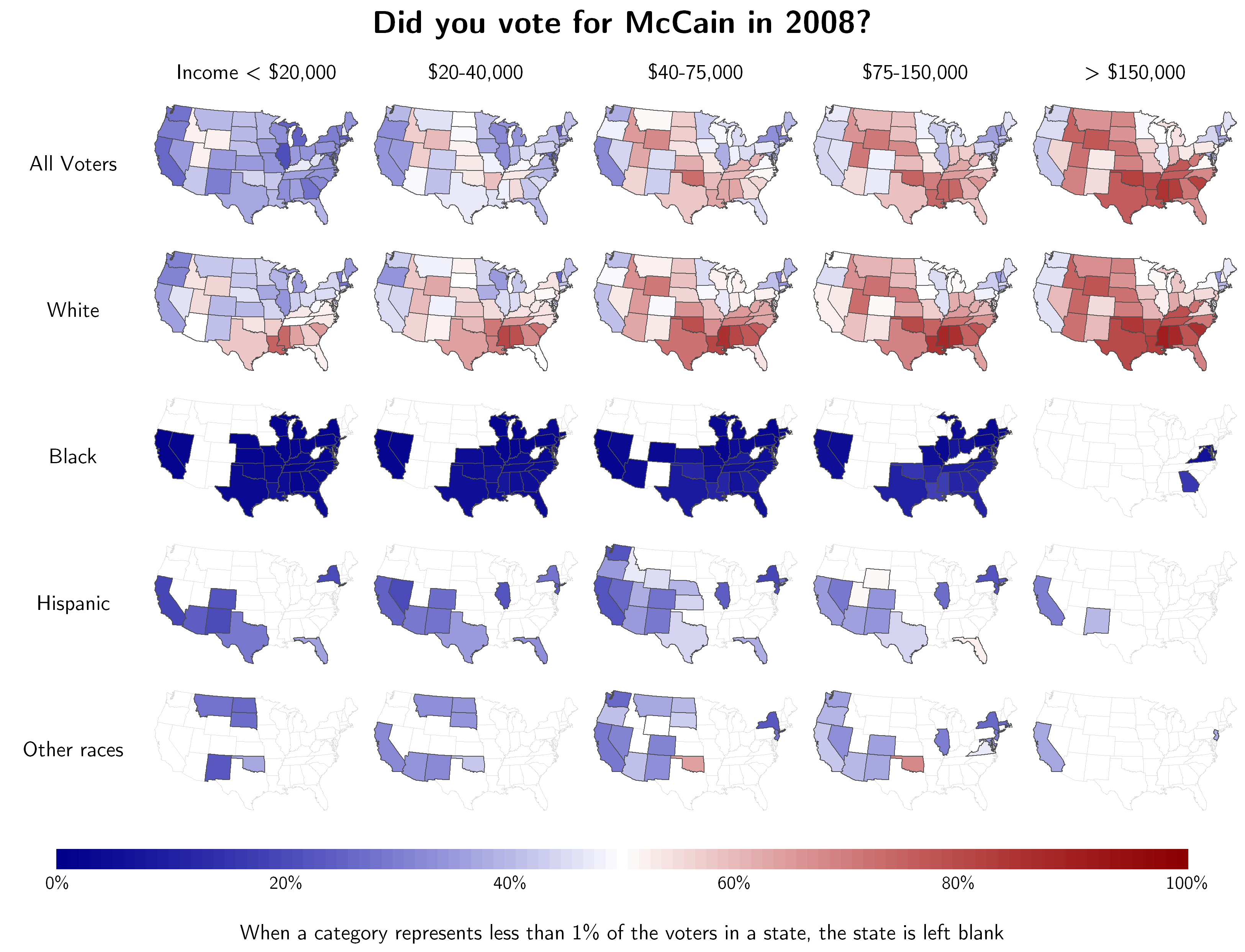}
  \end{center}
  \caption{\em Based on a model fitted to survey data: states won by John
    McCain and Barack Obama among different ethnic and income
    categories. States colored deep red and deep blue indicate clear McCain and
    Obama wins; pink and light blue represent wins by narrower margins, with a
    continuous range of shades going to gray for states estimated at exactly
    50/50.  The estimates shown here represent the culmination of months of
    effort, in which we fit increasingly complex models, at each stage checking
    the fit by comparing to data and then modifying aspects of the prior
    distribution and likelihood as appropriate.}
  \label{figure:1}
\end{figure}

\begin{figure}
  \resizebox{\textwidth}{!}{\includegraphics{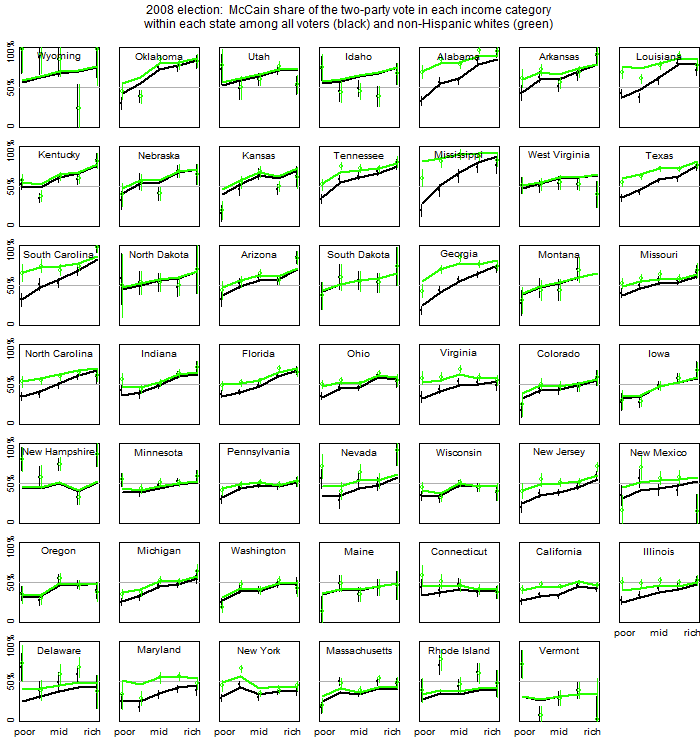}}
  \caption{\em Some of the data and fitted model used to make the maps shown in
    Figure \ref{figure:1}. Dots are weighted averages from pooled June-November
    Pew surveys; error bars show $\pm 1$ standard error bounds. Curves are
    estimated using multilevel models and have a standard error of about 3\% at
    each point. States are ordered in decreasing order of McCain vote (Alaska,
    Hawaii, and D.C. excluded). We fit a series of models to these data; only
    this last model fit the data well enough that we were satisfied. In working
    with larger datasets and studying more complex questions, we encounter
    increasing opportunities to check model fit and thus falsify in a way that
    is helpful for our research goals.}
  \label{figure:2}
\end{figure}

We demonstrate the hypothetico-deductive Bayesian modeling process with an
example from our recent applied research \citep{Gelman-Lee-Ghitza}. 
In recent years, American political
scientists have been increasingly interested in the connections between
politics and income inequality (see, e.g.,
\citealt{McCarty-Poole-Rosenthal-polarized}).  In our own contribution to this
literature, we estimated the attitudes of rich, middle-income, and poor voters
in each of the fifty states \citep{Gelman-et-al-red-state-blue-state}. As we
described in our article on the topic
\citep{rich-state-poor-state-red-state-blue-state}, we began by fitting a
varying-intercept logistic regression: modeling votes (coded as $y=1$ for votes
for the Republican presidential candidate or $y=0$ for Democratic votes) given
family income (coded in five categories from low to high as $x=-2, -1, 0, 1,
2$), using a model of the form $\Prob{y=1} = \mbox{logit}^{-1}(a_s + b x)$,
where $s$ indexes state of residence---the model is fit to survey
responses---and the varying intercepts $a_s$ correspond to some states being
more Republican-leaning than others. Thus, for example $a_s$ has a positive
value in a conservative state such as Utah and a negative value in a liberal
state such as California. The coefficient $b$ represents the ``slope'' of
income, and its positive value indicates that, within any state, richer voters
are more likely to vote Republican.

It turned out that this varying-intercept model did not fit our data, as we
learned by making graphs of the average survey response and fitted curves for
the different income categories within each state. We had to expand to a
varying-intercept, varying-slope model, $\Prob{y=1} = \mbox{logit}^{-1}(a_s +
b_s x)$, in which the slopes $b_s$ varied by state as well. This model
expansion led to a corresponding expansion in our understanding: we learned
that the gap in voting between rich and poor is much greater in poor states
such as Mississippi than in rich states such as Connecticut. Thus, the
polarization between rich and poor voters varied in important ways
geographically.

We found this not through any process of Bayesian induction but rather through
model checking. Bayesian inference was crucial, not for computing the posterior
probability that any particular model was true---we never actually did
that---but in allowing us to fit rich enough models in the first place that we
could study state-to-state variation, incorporating in our analysis relatively
small states such as Mississippi and Connecticut that did not have large
samples in our survey.  (\citealt{Gelman-Hill-data-analysis} review the
hierarchical models that allow such partial pooling.)

Life continues, though, and so do our statistical struggles. After the 2008
election, we wanted to make similar plots, but this time we found that even our
more complicated logistic regression model did not fit the data---especially
when we wanted to expand our model to estimate voting patterns for different
ethnic groups. Comparison of data to fit led to further model expansions,
leading to our current specification, which uses a varying-intercept,
varying-slope logistic regression as a baseline but allows for nonlinear and
even non-monotonic patterns on top of that. Figure \ref{figure:1} shows some of
our inferences in map form, while Figure \ref{figure:2} shows one of our diagnostics of data
and model fit.

The power of Bayesian inference here is {\em deductive}: given the data and some
model assumptions, it allows us to make lots of inferences, many of which can
be checked and potentially falsified. For example, look at New York state (in
the bottom row of Figure \ref{figure:2}): apparently, voters in the second
income category supported John McCain much more than did voters in neighboring
income groups in that state. This pattern is theoretically possible but it arouses
suspicion. A careful look at the graph reveals that this is a pattern in the
raw data which was moderated but not entirely smoothed away by our model. The
natural next step would be to examine data from other surveys. We may have
exhausted what we can learn from this particular dataset, and Bayesian
inference was a key tool in allowing us to do so.

\section{The Bayesian principal-agent problem}

Before returning to discussions of induction and falsification, we briefly
discuss some findings relating to Bayesian inference under misspecified models.
The key idea is that Bayesian inference for model selection---statements about
the posterior probabilities of candidate models---does not solve the problem of
learning from data about problems with existing models.

In economics, the ``principal-agent problem'' refers to the difficulty of
designing contracts or institutions which ensure that one selfish actor, the
``agent,'' will act in the interests of another, the ``principal,'' who cannot
monitor and sanction their agent without cost or error.  The problem is one of
aligning incentives, so that the agent serves itself by serving the principal
\citep{Eggertsson-economic-behavior-and-institutions}.  There is, as it were, a
Bayesian principal-agent problem as well.  The Bayesian agent is the
methodological fiction (now often approximated in software) of a creature with
a prior distribution over a well-defined hypothesis space $\Theta$, a
likelihood function $p(y|\theta)$, and conditioning as its sole mechanism of
learning and belief revision.  The principal is the actual statistician or
scientist.

The Bayesian agent's ideas are much more precise than the actual scientist's;
in particular, the Bayesian (in this formulation, with which we disagree) is
certain that {\em some} $\theta$ is the exact and complete truth, whereas the
scientist is not.\footnote{In claiming that ``the Bayesian'' is certain that
  some $\theta$ is the exact and complete truth, we are not claiming that
  actual Bayesian scientists or statisticians hold this view.  Rather, we are
  saying that this is implied by the philosophy we are attacking here.  All
  statisticians, Bayesian and otherwise, recognize that the philosophical
  position which ignores this approximation is problematic.}  At some point in
history, a statistician may well write down a model which he or she believes
contains all the systematic influences among properly-defined variables for the
system of interest, with correct functional forms and distributions of noise
terms.  This could happen, but we have never seen it, and in social science
we've never seen anything that comes close, either.  If nothing else, our own
experience suggests that however many different specifications we think of,
there are always others which had not occurred to us, but cannot be immediately
dismissed {\em a priori}, if only because they can be seen as alternative
approximations to the ones we made.  Yet the Bayesian agent is required to
start with a prior distribution whose support covers {\em all} alternatives
that could be considered.\footnote{It is also not at all clear that Savage and
  other founders of Bayesian decision theory ever thought that this principle
  should apply outside of the small worlds of artificially simplified and
  stylized problems---see \citet{Binmore-making-decisions-in-large-worlds}.
  But as scientists we care about the real, large world.}

This is not a small technical problem to be handled by adding a special value
of $\theta$, say $\theta^{\infty}$ standing for ``none of the above''; even if
one could calculate $p(y|\theta^{\infty})$, the likelihood of the data under
this catch-all hypothesis, this in general would {\em not} lead to just a small
correction to the posterior, but rather would have substantial effects
\citep{Fitelson-Thomason-implausible-theories}.  Fundamentally, the Bayesian
agent is limited by the fact that its beliefs always remain within the support
of its prior.  For the Bayesian agent, the truth must, so to speak, be always
already partially believed before it can become known.  This point is less than
clear in the usual treatments of Bayesian convergence, and so worth some
attention.

Classical results
\citep{Doob-bayesian-consistency,Schervish-theory-of-stats,Bayesian-consistency-for-stationary-models}
show that the Bayesian agent's posterior distribution will concentrate on the
truth with {\em prior} probability 1, provided some regularity conditions are
met.  Without diving into the measure-theoretic technicalities, the conditions
amount to (i) the truth is in the support of the prior, and (ii) the
information set is rich enough that some consistent estimator exists.  (See the
discussion in \citet[\S 7.4.1]{Schervish-theory-of-stats}.)  When the truth is
{\em not} in the support of the prior, the Bayesian agent still thinks that
Doob's theorem applies and assigns zero prior probability to the set of data
under which it does not converge on the truth.

The convergence behavior of Bayesian updating with a misspecified model can be
understood as follows
\citep{Berk-limiting-behavior-of-posterior,Berk-consistency,Kleijn-van-der-Vaart,CRS-dynamics-of-bayes}.
If the data are actually coming from a distribution $q$, then the
Kullback-Leibler divergence rate, or relative entropy rate, of the parameter
value $\theta$ is
\[
d(\theta) = \lim_{n\rightarrow\infty}{\frac{1}{n}\Expect{\log{\frac{p(y_1, y_2, \ldots y_n|\theta)}{q(y_1, y_2, \ldots y_n)}}}},
\]
with the expectation being taken under $q$.  (For details on when the limit
exists, see \citealt{Gray-entropy}.)  Then, under not-too-onerous regularity
conditions, one can show \citep{CRS-dynamics-of-bayes} that
\[
p(\theta|y_1, y_2, \ldots y_n) \approx p(\theta) \exp{\left\{-n(d(\theta) - d^*)\right\}},
\]
with $d^*$ being the essential infimum of the divergence rate.  More exactly,
\[
-\frac{1}{n}\log{p(\theta|y_1, y_2, \ldots y_n)} \rightarrow d(\theta) - d^*,
\]
$q$-almost-surely.  Thus the posterior distribution comes to concentrate on the
parts of the prior support which have the lowest values of $d(\theta)$ and the
highest expected likelihood.\footnote{More precisely, regions of $\Theta$ where
  $d(\theta) > d^*$ tend to have exponentially small posterior probability;
  this statement covers situations like $d(\theta)$ only approaching its
  essential infimum as $\|\theta\| \rightarrow \infty$, etc.  See
  \citet{CRS-dynamics-of-bayes} for details.}  There is a geometric sense in
which these parts of the parameter space are closest approaches to the truth
within the support of the prior \citep{Kass-Vos}, but they may or may not be
close to the truth in the sense of giving accurate values for parameters of
scientific interest.  They may not even be the parameter values which give the
best predictions
\citep{Grunwald-Langford-suboptimal-bayes,Muller-risk-of-Bayesian-inference}.
In fact, one cannot even guarantee that the posterior will concentrate on a
single value of $\theta$ at all; if $d(\theta)$ has multiple global minima, the
posterior can alternate between (concentrating around) them forever
\citep{Berk-limiting-behavior-of-posterior}.

To sum up, what Bayesian updating does when the model is false (that is, in
reality, always) is to try to concentrate the posterior on the best attainable
approximations to the distribution of the data, ``best'' being measured by
likelihood.  But depending on {\em how} the model is misspecified, and how
$\theta$ represents the parameters of scientific interest, the impact of
misspecification on inferring the latter can range from non-existent to
profound.\footnote{\citet{White-specification-analysis} gives examples of
  econometric models where the influence of mis-specification on the parameters
  of interest runs through this whole range, though only considering maximum
  likelihood and maximum quasi-likelihood estimation.}  Since we are quite
sure our models are wrong, we need to check whether the misspecification is so
bad that inferences regarding the scientific parameters are in trouble.  It is
by this non-Bayesian checking of Bayesian models that we solve our
principal-agent problem.

\section{Model checking}
\label{sec:model-checking}

In our view, a key part of Bayesian data analysis is model checking, which is
where there are links to falsificationism.  In particular, we emphasize the
role of posterior predictive checks, creating simulations and comparing the
simulated and actual data.  Again, we are following the lead of \citet{Box-sampling-and-Bayes},
\citet{Rubin-bayesianly-justifiable} and others, also mixing in a bit of
\citet{Tukey-EDA} in that we generally focus on visual comparisons \citep[ch.\
6]{Gelman-et-al-Bayesian-data-analysis}.

Here's how this works.  A Bayesian model gives us a joint distribution for the
parameters $\theta$ and the observables $y$.  This implies a marginal
distribution for the data,
\[
p(y) = \int{p(y|\theta) p(\theta) d\theta}.
\]
If we have observed data $y$, the prior distribution $p(\theta)$ shifts to the
posterior distribution $p(\theta|y)$, and so a different distribution of
observables,
\[
p(y^{\rm rep}|y) = \int{p(y^{\rm rep}|\theta) p(\theta|y) d\theta},
\]
where we use the $y^{\rm rep}$ to indicate hypothetical alternative or future
data, a replicated data set of the same size and shape as the original $y$,
generated under the assumption that the fitted model, prior and likelihood
both, is true.  By simulating from the posterior distribution of $y^{\rm rep}$,
we see what typical realizations of the fitted model are like, and in
particular whether the observed dataset is the kind of thing that the fitted model
produces with reasonably high probability.\footnote{For notational simplicity, we leave
  out the possibility of generating new values of the hidden variables
  $\tilde{y}$ and set aside choices of which parameters to vary and which to
  hold fixed in the replications; see \citet{Gelman-Meng-Stern}.} 

If we summarize the data with a test statistic $T(y)$, we can perform graphical comparisons with replicated data. In practice, we recommend graphical comparisons (as illustrated by our example above), but for continuity with much of the statistical literature, we focus here on
$p$-values,
\[
\Prob{T(y^{\rm rep}) > T(y)| y},
\]
which can be approximated to arbitrary accuracy as soon as we can simulate
$y^{\rm rep}$.  (This is a valid posterior probability in the model, and its
interpretation is no more problematic than that of any other probability in a
Bayesian model.)  In practice, we find graphical test summaries more
illuminating than $p$-values, but in considering ideas of (probabilistic)
falsification, it can be helpful to think about numerical test
statistics.\footnote{There is some controversy in the literature about whether
  posterior predictive checks have too little power to be useful statistical
  tools
  \citep{Bayarri-Berger-p-values-for-composite-null-models,Bayarri-Berger-interplay},
  how they might be modified to increase their power
  \citep{Robins-van-der-Vaart-VV,Fraser-Rousseau-p-values}, whether some form
  of empirical prior predictive check might not be better
  \citep{Bayarri-Castellanos-checking}, etc.  This is not the place to re-hash
  this debate over the interpretation or calculation of various Bayesian
  tail-area probabilities \citep{Gelman-on-Bayarri-Castellanos}.  Rather, the
  salient fact is that all participants in the debate agree on {\em why} the
  tail-area probabilities are relevant: they make it possible to reject a
  Bayesian model without recourse to a specific alternative.  All participants
  thus {\em disagree} with the standard inductive view, which reduces inference
  to the probability that a hypothesis is true, and are simply trying to find
  the most convenient and informative way to check Bayesian models.}

Under the usual understanding that $T$ is chosen so large values indicate poor
fits, these $p$-values work rather like classical ones
\citep{Mayo-error,Mayo-Cox-frequentist}---they in fact are generalizations of
classical $p$-values, merely replacing point estimates of parameters $\theta$
with averages over the posterior distribution---and their basic logic is one of
falsification.  A very low $p$-value says that it is very improbable, under the
model, to get data as extreme along the $T$-dimension as the actual $y$; we are
seeing something which would be very improbable if the model were true.  On the
other hand a high $p$-value merely indicates that $T(y)$ is an aspect of the
data which would be unsurprising if the model is true.  Whether this is
evidence {\em for} the usefulness of the model depends how likely it is to get
such a high $p$-value when the model is false: the ``severity'' of the test, in
the terminology of \citet{Mayo-error} and \citet{Mayo-Cox-frequentist}.

Put a little more abstractly, the hypothesized model makes certain
probabilistic assumptions, from which other probabilistic implications follow
deductively.  Simulation works out what those implications are, and tests check
whether the data conform to them.  Extreme $p$-values indicate that the data
violate regularities implied by the model, or approach doing so.  If these were
strict violations of deterministic implications, we could just apply {\em modus
  tollens} to conclude that the model was wrong; as it is, we nonetheless have
evidence and probabilities.  Our view of model checking, then, is firmly in the
long hypothetico-deductive tradition, running from \citet{Popper-logic} back
through \citet{Bernard} and beyond \citep{Laudan-science-and-hypothesis}.  A
more direct influence on our thinking about these matters is the work of
\citet{Jaynes-book}, who illustrated how we may learn the most when we find
that our model does not fit the data---that is, when it is falsified---because
then we have found a problem with our model's assumptions.\footnote{A similar
  point was expressed by the sociologist and social historian Charles
  \citet[p.\ 597]{Tilly-formal-representations}, writing from a very different
  disciplinary background: ``Most social researchers learn more from being
  wrong than from being right---provided they then recognize that they were
  wrong, see why they were wrong, and go on to improve their arguments.  Post
  hoc interpretation of data minimizes the opportunity to recognize
  contradictions between arguments and evidence, while adoption of formalisms
  increases that opportunity.  Formalisms blindly followed induce blindness.
  Intelligently adopted, however, they improve vision.  Being obliged to spell
  out the argument, check its logical implications, and examine whether the
  evidence conforms to the argument promotes both visual acuity and
  intellectual responsibility.''} And the better our probability model encodes
our {\em scientific} or {\em substantive} assumptions, the more we learn from
specific falsification.

In this connection, the prior distribution $p(\theta)$
is one of the assumptions of the model and does not need to represent the
statistician's personal degree of belief in alternative parameter values.  The
prior is connected to the data, and so is potentially testable, via the
posterior predictive distribution of future data $y^{\rm rep}$:
\begin{eqnarray*}
p(y^{\rm rep}|y) & = & \int{p(y^{\rm rep}|\theta) p(\theta|y) d\theta}\\
& = & \int{p(y^{\rm rep}|\theta) \frac{p(y|\theta)p(\theta)}{\int\!{p(y|\theta^{\prime})p(\theta^{\prime})d\theta^{\prime}}}d\theta}.
\end{eqnarray*}
The prior distribution thus has implications for the distribution of replicated
data, and so can be checked using the type of tests we have described and
illustrated above.\footnote{Admittedly, the prior only has observable
  implications in conjunction with the likelihood, but for a Bayesian the
  reverse is also true.} When it makes sense to think of further data coming
from the same source, as in certain kinds of sampling, time-series or
longitudinal problems, the prior also has implications for these new data
(through the same formula as above, changing the interpretation of $y^{\rm
  rep}$), and so becomes testable in a second way.  There is thus a connection
between the model-checking aspect of Bayesian data analysis and
``prequentialism'' \citep{Dawid-Vovk-prequential,Grunwald-on-MDL}, but
exploring that would take us too far afield.

One advantage of recognizing that the prior distribution is a testable part of
a Bayesian model is that it clarifies the role of the prior in inference, and
where it comes from.  To reiterate, it is hard to claim that the prior
distributions used in applied work represent statisticians' states of knowledge
and belief before examining their data, if only because most statisticians do
not believe their models are true, so their prior degree of belief in all of
$\Theta$ is not 1 but 0.  The prior distribution is more like a regularization
device, akin to the penalization terms added to the sum of squared errors when
doing ridge regression and the lasso \citep{tEoSL-2nd} or spline smoothing
\citep{Wahba-spline-models}.  All such devices exploit a sensitivity-stability
tradeoff: they stabilize estimates and predictions by making fitted models less
sensitive to certain details of the data.  Using an informative prior
distribution (even if only weakly informative, as in
\citet{Gelman-et-al-default-prior-for-logistic}) makes our estimates less
sensitive to the data than, say, maximum-likelihood estimates would be, which
can be a net gain.

Because we see the prior distribution as a testable part of the Bayesian model,
we do not need to follow Jaynes in trying to devise unique, objectively-correct
prior distribution for each situation---an enterprise with an uninspiring
track record \citep{Kass-Wasserman-selection-of-priors}, even leaving aside
doubts about Jaynes's specific proposal
\citep{Seidenfeld-why-i-am-not,Seidenfeld-entropy-and-uncertainty,Csiszar-on-maxent,Uffink-max-ent-and-consistency,Uffink-constraint-rule}.
To put it even more succinctly, ``the model,'' for a Bayesian, is the
combination of the prior distribution and the likelihood, each of which
represents some compromise among scientific knowledge, mathematical
convenience, and computational tractability.

This gives us a lot of flexibility in modeling.  We do not have to worry about
making our prior distributions match our subjective beliefs, still less about
our model containing all possible truths. Instead we make some assumptions,
state them clearly, see what they imply, and check the implications.  This
applies just much to the prior distribution as it does to the parts of the
model showing up in the likelihood function.

\subsection{Testing to reveal problems with a model}

We are not interested in falsifying our model for its own sake---among other
things, having built it ourselves, we know all the shortcuts taken in doing so,
and can already be morally certain it is false.  With enough data, we can
certainly detect departures from the model---this is why, for example, statistical
folklore says that the chi-squared statistic is ultimately a measure of sample
size (cf.\ \citealt{Lindsay-Liu-model-assessment}).  As writers such as
\citet[ch.\ 3]{Giere-explaining} explain, the hypothesis linking mathematical
models to empirical data is not that the data-generating process is exactly
isomorphic to the model, but that the data source resembles the model closely
enough, in the respects which matter to us, that reasoning based on the model
will be reliable.  Such reliability does not require complete fidelity to the
model.

The goal of model checking, then, is not to demonstrate the foregone conclusion
of falsity as such, but rather to learn how, in particular, this model fails
\citep{Gelman-goodness-of-fit}.\footnote{In addition, no model is safe from
  criticism, even if it ``passes'' all possible checks.  Modern Bayesian models
  in particular are full of unobserved, latent, and unobservable variables, and
  nonidentifiability is an inevitable concern in assessing such models; see,
  for example,
  \citet{Gustafson-on-model-expansion,Vansteelandt-et-al-on-uncertainty,Greenland-nonidentified-bias-sources}.
  We find it somewhat dubious to claim that simply putting a prior distribution
  on non-identified quantities somehow resolves the problem; the ``bounds'' or
  ``partial identification'' approach, pioneered by
  \citet{Manski-identification-for-prediction}, seems to be in better accord
  with scientific norms of explicitly acknowledging uncertainty (see also
  \citealt{Vansteelandt-et-al-on-uncertainty,Greenland-nonidentified-bias-sources}).}
When we find such particular failures, they tell us how the model must be
improved; when severe tests cannot find them, the inferences we draw about
those aspects of the real world from our fitted model become more credible.  In
designing a {\em good} test for model checking, we are interested in finding
particular errors which, if present, would mess up particular inferences, and
devise a test statistic which is sensitive to this sort of mis-specification.
This process of examining, and ruling out, possible errors or
mis-specifications is of course very much in line with the ``eliminative
induction'' advocated by \citet[ch.\ 7]{Kitcher-advancement}\footnote{Despite
  the name, this is, as Kitcher notes, actually a deductive argument.}

All models will have errors of approximation.  Statistical models, however,
typically assert that their errors of approximation will be unsystematic and
patternless---``noise'' \citep{Spanos-Ptolmey-vs-Kepler}.  Testing this can be
valuable in revising the model.  In looking at the red-state/blue-state
example, for instance, we concluded that the varying slopes mattered not just
because of the magnitudes of departures from the equal-slope assumption, but
also because there was a pattern, with richer states tending to have shallower
slopes.

What we are advocating, then, is what \citet{Cox-Hinkley} call ``pure
significance testing,'' in which certain of the model's implications are
compared directly to the data, rather than entering into a contest with some
alternative model.  This is, we think, more in line with what actually happens
in science, where it can become clear that even large-scale theories are
in serious trouble and cannot be accepted unmodified even if there is no
alternative available yet.  A classical instance is the status of Newtonian
physics at the beginning of the 20th century, where there were enough
difficulties---the Michaelson-Morley effect, anomalies in the orbit of Mercury,
the photoelectric effect, the black-body paradox, the stability of charged
matter, etc.---that it was clear, even before relativity and quantum mechanics,
that something would have to give.  Even today, our current best theories of
fundamental physics, namely general relativity and the standard model of
particle physics, an instance of quantum field theory, are universally agreed
to be ultimately wrong, not least because they are mutually incompatible, and
recognizing this does not require that one have a replacement theory
\citep{Weinberg-was-is-qft}.

\subsection{Connection to non-Bayesian model checking}
\label{sec:connection-to-non-bayesian-checking}

Many of these ideas about model checking are not unique to Bayesian data
analysis and are used more or less explicitly by many communities of
practitioners working with complex stochastic models
(\citealt{Ripley-spatial-processes,Guttorp-scientific}).  The reasoning is the
same: a model is a story of how the data could have been generated;
the fitted model should therefore be able to generate synthetic data that look
like the real data; failures to do so in important ways indicate faults in the
model.

For instance, simulation-based model checking is now widely accepted for
assessing the goodness of fit of statistical models of social networks
\citep{Hunter-Goodreau-Handcock-gof-of-social-networks}.  That community was
pushed toward predictive model checking by the observation that many model
specifications were ``degenerate'' in various ways
\citep{Handcock-assessing-degeneracy}.  For example, under certain
exponential-family network models, the maximum likelihood estimate gave a
distribution over networks which was bimodal, with both modes being very
different from observed networks, but located so that the expected value of the
sufficient statistics matched observations.  It was thus clear that these
specifications could not be right even before more adequate specifications were
developed \citep{Snijders-Pattison-Robins-Handcock-new-specification}.

At a more philosophical level, the idea that a central task of statistical
analysis is the search for specific, consequential errors has been forcefully
advocated by
\citet{Mayo-error}, \citet{Mayo-Cox-frequentist,Mayo-Spanos-methodology-and-practice}, and \citet{Mayo-Spanos-post-data}.
Mayo has placed a special emphasis on the idea of {\em severe} testing---a
model being severely tested if it passes a probe which had a high probability
of detecting an error if it is present.  (The exact definition of a test's
severity is related to, but not quite, that of its power; see
\citealt{Mayo-error} or \citealt{Mayo-Spanos-post-data} for extensive
discussions.)  Something like this is implicit in discussions about
the relative merits of particular posterior predictive checks (which can also
be framed non-Bayesianly as graphical hypothesis tests based on the parametric
bootstrap).

Our contribution here is to connect this hypothetico-deductive philosophy to
Bayesian data analysis, going beyond the evaluation of Bayesian methods based
on their frequency properties (as recommended by
\citet{Rubin-bayesianly-justifiable}, \citet{Wasserman-in-Bayesian-Analysis},
among others) to emphasize the learning that comes from the discovery of
systematic differences between model and data.  At the very least, we hope this
paper will motivate philosophers of hypothetico-deductive inference to take a
more serious look at Bayesian data analysis (as distinct from Bayesian theory)
and, conversely, to motivate philosophically-minded Bayesian statisticians to
consider alternatives to the inductive interpretation of Bayesian learning.

\subsection{Why not just compare the posterior probabilities of different models?}

As mentioned above, the standard view of scientific learning in the Bayesian
community is, roughly, that posterior odds of the models under
consideration are compared, given the current data.\footnote{Some would prefer to compare
  the modification of those odds called the Bayes factor
  \citep{Kass-Raftery-bayes-factors}.  Everything we have to say about
  posterior odds carries over to Bayes factors with few changes.}  When
Bayesian data analysis is understood as simply getting the posterior
distribution, it is held that ``pure significance tests have no role to play in
the Bayesian framework'' \citep[p.\ 218]{Schervish-theory-of-stats}.  The
dismissal rests on the idea that the prior distribution can accurately reflect
our actual knowledge and beliefs.\footnote{As \citet{Schervish-theory-of-stats}
  continues: ``If the [parameter space $\Theta$] describes all of the
  probability distributions one is willing to entertain, then one cannot reject
  [$\Theta$] without rejecting probability models altogether.  If one is
  willing to entertain models not in [$\Theta$], then one needs to take them
  into account'' by enlarging $\Theta$, and computing the posterior
  distribution over the enlarged space.}  At the risk of boring the reader by
repetition, there is just no way we can ever have any hope of making $\Theta$
include all the probability distributions which might be correct, let alone
getting $p(\theta|y)$ if we did so, so this is deeply unhelpful advice.  The
main point where we disagree with many Bayesians is that we do not see Bayesian
methods as generally useful for giving the posterior probability that a model
is true, or the probability for preferring model A over model B, or
whatever.\footnote{There is a vast literature on Bayes factors, model
  comparison, model averaging, and the evaluation of posterior probabilities of
  models, and although we believe most of this work to be philosophically
  unsound (to the extent it is designed to be a direct vehicle for scientific
  learning), we recognize that these can be useful techniques.  Like all
  statistical methods, Bayesian and otherwise, these methods are summaries of
  available information that can be important data-analytic tools.  Even if
  none of a class of models is plausible as truth, and even if we aren't
  comfortable accepting posterior model probabilities as degrees of belief in
  alternative models, these probabilities can still be useful as tools for
  prediction and for understanding structure in data, as long as these
  probabilities are not taken too seriously.  See
  \citet{Raftery-bayesian-model-selection} for a discussion of the value of
  posterior model probabilities in social science research and
  \citet{Gelman-Rubin-contra-Raftery} for a discussion of their limitations,
  and \citet{Claeskens-Hjort-model-selection} for a general review of model
  selection.  (Some of the work on ``model-selection tests'' in econometrics
  (e.g.,
  \citealt{Vuong-testing-non-nested-hypotheses,Rivers-Vuong-model-selection-tests})
  is exempt from our strictures, as it tries to find which model is {\em
    closest} to the data-generating process, while allowing that all of the
  models may be mis-specified, but it would take us too far afield to discuss
  this work in detail.)} Beyond the philosophical difficulties, there are
technical problems with methods that purport to determine the posterior
probability of models, most notably that in models with continuous parameters,
aspects of the model that have essentially no effect on posterior inferences
{\em within} a model can have huge effects on the comparison of posterior
probability {\em among} models.\footnote{This problem has been called the
  Jeffreys-Lindley paradox and it is the subject of a large literature.
  Unfortunately (from our perspective) the problem has usually been studied by
  Bayesians with an eye toward ``solving'' it---that is, coming up with
  reasonable definitions that allow the computation of nondegenerate posterior
  probabilities for continuously-parameterized models---but we we think that
  this is really a problem without a solution; see \citet[sec.\
  6.7]{Gelman-et-al-Bayesian-data-analysis}.}  Bayesian inference is good for
deductive inference within a model we prefer to evaluate a model by comparing it to data.

In rehashing the well-known problems with computing Bayesian posterior
probabilities of models, we are not claiming that classical $p$-values are the
answer.  As is indicated by the literature on the Jeffreys-Lindley paradox
(notably \citealt{Berger-Sellke-testing-a-point-null}), $p$-values can
drastically overstate the evidence against a null hypothesis.  From our
model-building Bayesian perspective, the purpose of $p$-values (and model
checking more generally) is not to reject a null hypothesis but rather to
explore aspects of a model's misfit to data.

In practice, if we are in a setting where model A or model B might be true, we
are inclined not to do {\em model selection} among these specified options, or
even to perform {\em model averaging} over them (perhaps with a statement such
as, ``We assign 40\% of our posterior belief to A and 60\% to B'') but rather
to do {\em continuous model expansion} by forming a larger model that includes
both A and B as special cases.  For example,
\citet{Merrill-on-spatial-distribution-of-party-support} used electoral and
survey data from Norway and Sweden to compare two models of political ideology
and voting: the ``proximity model'' (in which you prefer the political party
that is closest to you in some space of issues and ideology) and the
``directional model'' (in which you like the parties that are in the same
direction as you in issue space, but with a stronger preference to parties
further from the center).  Rather than using the data to pick one model or the
other, we would prefer to think of a model in which voters consider both
proximity and directionality in forming their preferences
\citep{Gelman-on-Merrill}.

In the social sciences, it is rare for there to be an underlying theory that
can provide meaningful constraints on the functional form of the expected
relationships among variables, let alone the distribution of noise
terms.\footnote{See \citet{Manski-identification-for-prediction} for a critique
  of the econometric practice of making modeling assumptions (such as
  linearity) with no support in economic theory, simply to get
  identifiability.}  Taken to its limit, then, the idea of continuous model
expansion counsels social scientists to pretty much give up using parametric
statistical models in favor of nonparametric, infinite-dimensional models, advice which which the on-going rapid development of Bayesian nonparametrics
\citep{Ghosh-Ramamoorthi,Hjort-Holmes-Muller-Walker} makes increasingly
practical.  While we are certainly sympathetic to this, and believe a greater
use of nonparametric models in empirical research is desirable on its own
merits (cf.\ \citealp{Li-Racine-nonparametric-econometrics}), it is worth sounding a few notes of caution.

A technical, but important, point concerns the representation of uncertainty in
Bayesian nonparametrics.  In finite-dimensional problems, the use of the
posterior distribution to represent uncertainty is in part supported by the
Bernstein-von Mises phenomenon, which ensures that large-sample credible
regions are also confidence regions.  This simply fails in infinite-dimensional
situations
\citep{DDCox-on-Bayesian-nonparametric-regression,Freedman-on-bernstein-von-mises},
so that a naive use of the posterior distribution becomes unwise\footnote{Even
  in parametric problems, \citet{Muller-risk-of-Bayesian-inference} shows that
  mis-specification can lead credible intervals to have sub-optimal coverage
  properties, which however can be fixed by a modification to their usual
  calculation.}.  (Since we regard the prior and posterior distributions as
regularization devices, this is not especially troublesome for us.)  Relatedly,
the prior distribution in a Bayesian nonparametric model is a stochastic
process, always chosen for tractability
\citep{Ghosh-Ramamoorthi,Hjort-Holmes-Muller-Walker}, and any pretense of
representing an actual inquirer's beliefs abandoned.

Most fundamentally, switching to nonparametric models does not really resolve
the issue of needing to make approximations, and check their adequacy.  All
nonparametric models themselves embody assumptions such as conditional
independence which are hard to defend except as approximations.  Expanding our
prior distribution to embrace {\em all} the models which are actually
compatible with our prior knowledge would result in a mess we simply could not
work with, nor interpret if we could.  This being the case, we feel there is no
contradiction between our preference for continuous model expansion and our use
of {\em adequately checked} parametric models.\footnote{A different perspective
  --- common in econometrics (e.g.,
  \citealt{Wooldridge-cross-section-and-panel}) and machine learning (e.g.,
  \citealt{tEoSL-2nd}) --- reduces the importance of models of the data source,
  either by using robust procedures that are valid under departures from
  modeling assumptions, or by focusing on prediction and external validation.
  We recognise the theoretical and practical appeal of both these approaches,
  which can be relevant to Bayesian inference.  (For example,
  \citealt{Rubin-role-of-randomization} justifies random assignment from a
  Bayesian perspective as a tool for obtaining robust inferences.)  But it is
  not possible to work with {\em all} possible models when considering fully
  probabilistic methods --- that is, Bayesian inferences that are summarized by
  joint posterior distributions rather than point estimates or predictions.
  This difficulty may well be a motivation for shifting the foundations of
  statistics away from probability and scientific inference, and towards
  developing a technology of robust prediction.  (Even when prediction is the
  only goal, with limited data bias-variance considerations can make even
  mis-specified parametric models superior to non-parametric models.)  This
  however goes far beyond the scope of the present article, which aims merely
  to explicate the implicit philosophy guiding current practice.}

\subsection{Example: Estimating the effects of legislative redistricting}

We use one of our own experiences \citep{Gelman-King-redistricting} to
illustrate scientific progress through model rejection.  We began by fitting a
model comparing treated and control units---state legislatures, immediately
after redistricting or not---following the usual practice of assuming a
constant treatment effect (parallel regression lines in ``after''
vs. ``before'' plots, with the treatment effect representing the difference
between the lines).  In this example, the outcome was a measure of partisan
bias, with positive values representing state legislatures where the Democrats
were overrepresented (compared to how we estimated the Republicans would have
done with comparable vote shares) and negative values in states where the
Republicans were overrepresented.  A positive treatment effect here would
correspond to a redrawing of the district lines that favored the Democrats.

\begin{figure}
  \begin{center}
    \resizebox{.3\textwidth}{!}{\includegraphics{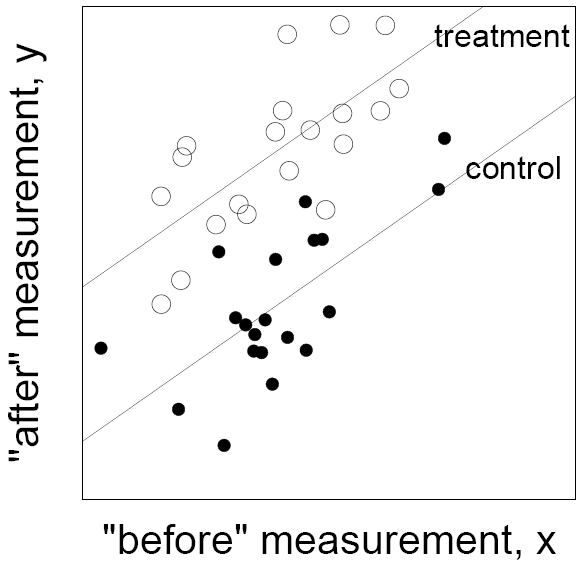}}
  \end{center}
  \caption{\em Sketch of the usual statistical model for before-after data.
    The difference between the fitted lines for the two groups is the estimated
    treatment effect.  The default is to regress the ``after'' measurement on
    the treatment indicator and the ``before'' measurement, thus implicitly
    assuming parallel lines.}
  \label{fig:beforeafterdefault}
\end{figure}

Figure \ref{fig:beforeafterdefault} shows the default model that we (and
others) typically use for estimating causal effects in before-after data.  We
fitted such a no-interaction model in our example too, but then we made some
graphs and realized that the model did not fit the data.  The line for the
control units actually had a much steeper slope than the treated units.  We fit
a new model, and it had a completely different story about what the treatment
effects meant.

\begin{figure}
  \begin{center}
    \resizebox{.6\textwidth}{!}{\includegraphics{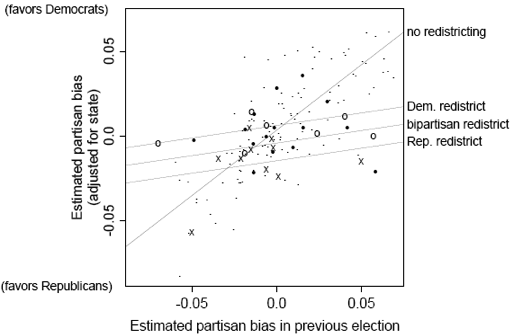}}
  \end{center}
\vspace{-.2in}
\caption{\em Effect of redistricting on partisan bias. Each symbol represents a
  state election year, with dots indicating controls (years with no
  redistricting) and the other symbols corresponding to different types of
  redistricting. As indicated by the fitted lines, the ``before'' value is much
  more predictive of the ``after'' value for the control cases than for the
  treated (redistricting) cases. The dominant effect of the treatment is to
  bring the expected value of partisan bias toward zero, and this effect would
  not be discovered with the usual approach (pictured in Figure
  \ref{fig:beforeafterdefault}, which is to fit a model assuming parallel
  regression lines for treated and control cases.}
  \label{figure:3}
\end{figure}

The graph for the new model with interactions is shown in Figure
\ref{figure:3}.  The largest effect of the treatment was not to benefit the
Democrats or Republicans (that is, to change the intercept in the regression,
shifting the fitted line up or down) but rather to change the slope of the line,
to reduce partisan bias.

Rejecting the constant-treatment-effect model and replacing by the interaction
model was, in retrospect, a crucial step in this research project.  This
pattern of higher before-after correlation in the control group than the
treated group is quite general \citep{Gelman-before-after-correlations}, but at
the time we did this study we discovered it only through the graph of model and data, which falsified the original model and motivated us to think of
something better.  In our experience, falsification is about plots and
predictive checks, not about Bayes factors or posterior probabilities of
candidate models.

The relevance of this example to the philosophy of statistics is that we began
by fitting the usual regression model with no interactions. Only after visually
checking the model fit---and thus falsifying it in a useful way without the
specification of any alternative---did we take the crucial next step of
including an interaction, which changed the whole direction of our research.
The shift was induced by a falsification---a bit of deductive inference from
the data and the earlier version of our model.  In this case the falsification came from a graph rather than a $p$-value, which in one way is just a technical issue, but in a larger sense is important in that the graph revealed not just a lack of fit but also a sense of the direction of the misfit, a refutation that sent us usefully in a direction of substantive model improvement.

\section{The question of induction}

As we mentioned at the beginning, Bayesian inference is often held to be
inductive in a way which classical statistics (following the Fisher or
Neyman-Pearson traditions) is not.  We need to address this, as we are arguing
that all these forms of statistical reasoning are better seen as
hypothetico-deductive.

The common core of various conceptions of induction is some form of inference
from particulars to the general---in the statistical context, presumably,
inference from the observations $y$ to parameters $\theta$ describing the
data-generating process.  But if {\em that} were all that was meant, then not
only is ``frequentist statistics a theory of inductive inference''
\citep{Mayo-Cox-frequentist}, but the whole range of guess-and-test behaviors
engaged in by animals \citet{HHNT-induction} are formalized in the
hypothetico-deductive method are also inductive.  Even the unpromising-sounding
procedure, ``Pick a model at random and keep it until its accumulated error
gets too big, then pick another model completely at random,'' would qualify
(and could work surprisingly well under some circumstances; cf.\
\citet{Ashby-design,Foster-Young-hypo-testing}).  So would utterly irrational
procedures (``pick a new random $\theta$ when the sum of the least significant
digits in $y$ is 13'').  Clearly something more is required, or at least
implied, by those claiming that Bayesian updating is inductive.

One possibility for that ``something more'' is to generalize the
truth-preserving property of valid deductive inferences: just as valid
deductions from true premises are themselves true, good inductions from true
observations should also be true, at least in the limit of increasing
evidence.\footnote{We owe this suggestion to conversation with Kevin Kelly;
  cf.\ \citet[esp.\ ch.\ 13]{Kelly-reliable}.}  This, however, is just the
requirement that our inferential procedures be consistent.  As discussed above,
using Bayes's rule is not sufficient to ensure consistency, nor is it
necessary.  In fact, every proof of Bayesian consistency known to us either
posits there is a consistent non-Bayesian procedure for the same problem, or
makes other assumptions which entail the existence of such a procedure.  In any
case, theorems establishing consistency of statistical procedures make {\em
  deductively valid} guarantees about these procedures---they are theorems,
after all---but do so on the basis of probabilistic assumptions linking future
events to past data.

It is also no good to say that what makes Bayesian updating inductive is its
conformity to some axiomatization of rationality.  If one accepts the
Kolmogorov axioms for probability, and the Savage axioms (or something like
them) for decision-making,\footnote{Despite his ideas on testing,
  \citet{Jaynes-book} was a prominent and emphatic advocate of the claim that
  Bayesian inference is the logic of inductive inference as such, but preferred
  to follow \citet{Cox-1946,Cox-algebra} rather than Savage.  See
  \citet{Halpern-Coxs-theorem-revisted} on the formal invalidity of Cox's
  proofs.} then updating by conditioning follows, and a prior belief state
$p(\theta)$ plus data $y$ {\em deductively} entail that the new belief state is
$p(\theta|y)$.  In any case, lots of learning procedures can be axiomatized
(all of them which can be implemented algorithmically, to start with).  To pick
{\em this} system, we would need to know that it produces good results (cf.\
\citealt{Manski-actualist}), and this returns us to previous problems.  To know that this
axioms systems leads us to approach the truth rather than become convinced of
falsehoods, for instance, is just the question of consistency again.

Karl Popper, the leading advocate of hypothetico-deductivism in the last
century, denied that induction was even possible; his attitude is
well-paraphrased by \citet{Greenland-on-Popper} as: ``we never use any argument
based on observed repetition of instances that does not also involve a
hypothesis that predicts both those repetitions and the unobserved instances of
interest.''  This is a recent instantiation of a tradition of anti-inductive
arguments that goes back to Hume, but also beyond him to
\citet{Ghazali-Tahafut} in the middle ages, and indeed to the ancient Skeptics
\citep{Kolakowski-positivism}.  As forcefully put by
\citet{Stove-popper-and-after,Stove-induction}, many apparent arguments against
this view of induction can be viewed as statements of abstract premises linking
both the observed data and unobserved instances---various versions of the
``uniformity of nature'' thesis have been popular, sometimes resolved into a
set of more detailed postulates, as in \citet[part VI,
ch. 9]{Russell-human-knowledge}, though Stove rather maliciously crafted a
parallel argument for the existence of ``angels, or something very much like
them.''\footnote{\citet{Stove-induction} further argues that induction by
  simple enumeration is reliable {\em without} making such assumptions, at
  least sometimes. However, his calculations make no sense unless his data are
  independent and identically distributed.}  As
\citet{Norton-material-induction} argues, these highly abstract premises are
both dubious and often superfluous for supporting the sort of actual inferences
scientists make --- ``inductions'' are supported not by their matching certain
formal criteria (as deductions are), but rather by material facts.  To
generalize about the melting point of bismuth (to use one of Norton's examples)
requires very few samples, provided we accept certain facts about the
homogeneity of the physical properties of elemental substances; whether nature
in general is uniform is not really at issue\footnote{ Within environments
  where such premises hold, it may of course be adaptive for organisms to
  develop inductive propensities, whose scope would be more or less tied to the
  domain of the relevant material premises.  \citet{Adapted-mind} develops this
  theme with reference to the evolution of domain-specific mechanisms of
  learning and induction;
  \citet{Gigerenzer-adaptive-thinking,Gigerenzer-Todd-heuristics} consider
  proximate mechanisms and ecological aspects, and \citet{HHNT-induction}
  proposes a unified framework for modeling such inductive propensities in
  terms of generate-and-test processes.  All of this, however, is more within
  the field of psychology than either statistics or philosophy, as (to
  paraphrase the philosopher Ian \citealt{Hacking-intro-to-prob}) it does not
  so much solve the problem of induction as evade it.}.

Simply put, we think the anti-inductivist view is pretty much right, but that
statistical models are tools that let us draw inductive inferences on a
deductive background.  Most directly, random sampling allows us to learn about
unsampled people (unobserved balls in an urn, as it were), but such inference,
however inductive it may appear, relies not any axiom of induction but rather
on deductions from the statistical properties of random samples, and the
ability to actually conduct such sampling.  The appropriate design
depends on many contingent material facts about the system we are studying,
exactly as Norton argues.

Some results in statistical learning theory establish that certain procedures
are ``probably approximately correct'' in what's called a ``distribution-free''
manner
\citep{Bousquet-Boucheron-Lugosi,Vidyasagar-on-learning-and-generalization};
some of these results embrace Bayesian updating \citep{McAllister-PAC-Bayes}.
But, here, ``distribution free'' just means ``holding uniformly over all
distributions in a very large class,'' for example requiring the data to be
independent and identically distributed, or from a stationary, mixing
stochastic process.  Another branch of learning theory does avoid making any
probabilistic assumptions, getting results which hold universally across all
possible data sets, and again these results apply to Bayesian updating, at
least over some parameter spaces \citep{prediction-learning-and-games}.
However, these results are all of the form ``in retrospect, the posterior
predictive distribution will have predicted almost as well as the best
individual model could have done,'' speaking entirely about performance on the
past training data and revealing nothing about extrapolation to so-far
unobserved cases.

To sum up, one is free to describe statistical inference as a theory of
inductive logic, but these would be inductions which are deductively guaranteed
by the probabilistic assumptions of stochastic models.  We can see no
interesting and correct sense in which Bayesian statistics is a logic of
induction which does not equally imply that frequentist statistics is also a
theory of inductive inference (cf. \citealp{Mayo-Cox-frequentist}), which is to
say, not very inductive at all.

\section{What About Popper and Kuhn?}

The two most famous modern philosophers of science are undoubtedly Karl
\citet{Popper-logic} and Thomas \citet{Kuhn-structure}, and if statisticians
(like other non-philosophers) know about philosophy of science at all, it is
generally some version of their ideas.  It may therefore help readers for see
how our ideas relate to theirs.  We do not pretend that our sketch fully
portrays these figures, let alone the literatures of exegesis and controversy
they inspired, or even how the philosophy of science has moved on since 1970.

Popper's key idea was that of ``falsification,'' or ``conjectures and
refutations.''  The inspiring example, for Popper, was the replacement of
classical physics, after several centuries as the core of the best-established
science, by modern physics, especially the replacement of Newtonian gravitation
by Einstein's general relativity.  Science, for Popper, advances by scientists
advancing theories which make strong, wide-ranging predictions capable of being
refuted by observations.  A good experiment or observational study is one which
tests a specific theory (or theories) by confronting their predictions with
data in such a way that a match is not automatically assured; good studies are
designed with theories in mind, to give them a chance to fail.  Theories which
conflict with any evidence must be rejected, since a single counter-example
implies that a generalization is false.  Theories which are not falsifiable by
any conceivable evidence are, for Popper, simply not scientific, though they
may have other virtues.\footnote{This ``demarcation criterion'' has received a
  lot of criticism, much of it justified.  The question of what makes something
  ``scientific'' is fortunately not one we have to answer; cf.\ \citet[chs.\
  11--12]{Laudan-beyond} and \citet{Ziman-real-science}.}  Even those
falsifiable theories which have survived contact with data so far must be
regarded as more or less provisional, since no finite amount of data can ever
establish a generalization, nor is there any non-circular principle of
induction which could let us regard theories which are compatible with lots of
evidence are probably true.\footnote{Popper tried to work out notions of
  ``corroboration'' and increasing truth content, or ``verisimilitude,'' that
  fit with these stances, but these are generally regarded as failures.}  Since
people are fallible, and often obstinate and overly fond of their own ideas,
the objectivity of the process which tests conjectures lies not in the
emotional detachment and impartiality of individual scientists, but rather in
the scientific community being organized in certain ways, with certain
institutions, norms and traditions, so that individuals' prejudices more or
less wash out \citep[chs.\ 23--24]{Popper-open-society}.

Clearly, we find much here to agree with, especially the general
hypothetico-deductive view of scientific method and the anti-inductivist
stance.  On the other hand, Popper's specific ideas about testing require, at
the least, substantial modification.  His idea of a test comes down to the rule
of deduction which says that if $p$ implies $q$, and $q$ is false, then $p$
must be false, with the roles of $p$ and $q$ being played by hypotheses and
data, respectively.  This is plainly inadequate for statistical hypotheses,
yet, as critics have noted since \citet{Braithwaite-scientific-explanation} at
least, he oddly ignored the theory of statistical hypothesis
testing.\footnote{We have generally found Popper's ideas on probability and
  statistics to be of little use and will not discuss them here.}  It is possible to do better, both through standard hypothesis tests and the
kind of predictive checks we have described.  In particular, as
\citet{Mayo-error} has emphasized, it is vital to consider the {\em severity}
of tests, their capacity to detect violations of hypotheses when they are
present.

Popper tried to say how science {\em ought} to work, supplemented by arguments
that his ideals could at least be approximated and often had been.  Kuhn's
work, by contrast, account, in contrast, was much more an attempt to describe
how science had, in point of historical fact, developed, supported by arguments
that alternatives were infeasible, from which some morals might be drawn.  His
central idea was that of a ``paradigm,'' a scientific problem and
its solution which served as a model or exemplar, so that solutions to other
problems could be developed in imitation of it.\footnote{Examples include
  Newton's deduction of Kepler's laws of planetary motion and other facts of
  astronomy from the inverse square law of gravitation, or Planck's derivation
  of the black-body radiation distribution from Boltzmann's statistical
  mechanics and the quantization of the electromagnetic field.  An internal
  example for statistics might be the way the Neyman-Pearson lemma inspired the
  search for uniformly most powerful tests in a variety of complicated situations.}
Paradigms come along with presuppositions about the terms available for
describing problems and their solutions, what counts as a valid problem, what
counts as a solution, background assumptions which can be taken as a matter of
course, etc.  Once a scientific community accepts a paradigm and all that goes
with it, its members can communicate with one another, and get on with the
business of ``puzzle solving,'' rather than arguing about what they should be
doing.  Such ``normal science'' includes a certain amount of developing and testing of hypotheses but leaves the central presuppositions of the paradigm
unquestioned.

During periods of normal science, according to Kuhn, there will always be some
``ano\-ma\-lies''---things within the domain of the paradigm which it currently
cannot explain, or even seem to refute its assumptions.  These are generally
ignored, or at most regarded as problems which somebody ought to investigate
eventually.  (Is a special adjustment for odd local circumstances called for?
Might there be some clever calculational trick which fixes things?  How sound
are those anomalous observations?)  More formally, Kuhn invokes the
``Quine-Duhem thesis'' \citep{Quine-logical,Duhem-aim-and-structure}. A
paradigm only makes predictions about observations in conjunction with
``auxiliary'' hypotheses about specific circumstances, measurement procedures,
etc.  If the predictions are wrong, Quine and Duhem claimed that one is always
free to fix the blame on the auxiliary hypotheses, and preserve belief in the
core assumptions of the paradigm ``come what may.''\footnote{This thesis can be
  attacked from many directions, perhaps the most vulnerable being that one can
  often find multiple lines of evidence which bear on either the main
  principles or the auxiliary hypotheses {\em separately}, thereby localizing
  the problems
  \citep{Glymour-theory-and-evidence,Kitcher-advancement,Laudan-beyond,Mayo-error}.}
The Quine-Duhem thesis was also used by \citet{Lakatos-papers} as part of his
``methodology of scientific research programmes,'' a falsificationism more
historically oriented than Popper's distinguishing between progressive
development of auxiliary hypotheses and degenerate research programs where
auxiliaries become {\em ad hoc} devices for saving core assumptions from data.

According to Kuhn, however, anomalies can accumulate, becoming so serious as to
create a crisis for the paradigm, beginning a period of ``revolutionary
science.''  It is then that a new paradigm can form, one which is generally
``incommensurable'' with the old: it makes different presuppositions, takes a
different problem and its solution as exemplars, re-defines the meaning of
terms.  Kuhn insisted that scientists who retain the old paradigm are not being
irrational, because (by Quine-Duhem) they can always explain away the anomalies
{\em somehow}; but neither are the scientists who embrace and develop the new
paradigm being irrational.  Switching to the new paradigm is more like a
bistable illusion flipping (the apparent duck becomes an obvious rabbit) than
any process of ratiocination governed by sound rules of
method.\footnote{\citet{Salmon-Kuhn-meets-Bayes} proposed a connection between
  Kuhn and Bayesian reasoning, suggesting that the choice between paradigms
  could be made rationally by using Bayes's rule to compute their posterior
  probabilities, with the prior probabilities for the paradigms encoding such
  things as preferences for parsimony.  This has at least three big problems.
  First, all our earlier objections to using posterior probabilities to chose
  between theories apply, with all the more force because every paradigm is
  compatible with a broad range of specific theories.  Second, devising priors
  encoding those methodological preferences---particularly a non-vacuous
  preference for parsimony---is hard to impossible in practice
  \citep{Kelly-simplicity-truth-probability}.  Third, it implies a truly
  remarkable form of Platonism: for scientists to give a paradigm positive
  posterior probability, they must, by Bayes's rule, have always given it
  strictly positive prior probability, {\em even before having encountered a
    statement of the paradigm}.\label{note:salmon}}


In some way, Kuhn's distinction between normal and revolutionary science is
analogous to the distinction between learning within a Bayesian model, and
checking the model as preparation to discard or expand it.  Just as the work of normal science proceeds within the
presuppositions of the paradigm, updating a posterior distribution by
conditioning on new data takes the assumptions embodied in the prior
distribution and the likelihood function as unchallengeable truths.  Model
checking, on the other hand, corresponds to the identification of anomalies,
with a switch to a new model when they become intolerable.  Even the problems
with translations between paradigms have something of a counterpart in
statistical practice; for example, the intercept coefficients in a
varying-intercept, constant-slope regression model have a somewhat different
meaning than do the intercepts in a varying-slope model.  We do not want to
push the analogy too far, however, since most model checking and model
re-formulation would by Kuhn have been regarded as puzzle-solving within a
single paradigm, and his views of how people switch between paradigms are, as
we just saw, rather different.

Kuhn's ideas about scientific revolutions are famous because they raise so many
disturbing questions about the scientific enterprise.  For instance, there has
been considerable controversy over whether Kuhn believed in any notion of
scientific progress, and over whether or not he should have, given his theory.
Yet detailed historical case studies \citep{scrutinizing-science} have shown
that Kuhn's picture of sharp breaks between normal and revolutionary science is
hard to sustain.  (Arguably this is true even of \citealp{Kuhn-copernican}.)
The leads to a tendency, already remarked by \citet[pp.\
112--17]{Toulmin-human-understanding}, to either expand paradigms or to shrink
them.  Expanding paradigms into persistent and all-embracing, because 
abstract and vague, bodies of ideas lets one preserve the idea of abrupt breaks
in thought, but makes them rare and leaves almost everything to
puzzle-solving normal science.  (In the limit, there has only been one paradigm
in astronomy since the Mesopotamians, something like ``Many lights in the night
sky are objects which are very large but very far away, and they move in
interrelated, mathematically-describable, discernible patterns.'')  This
corresponds, we might say, to relentlessly enlarging the support of the prior.
The other alternative is to shrink paradigms into increasingly concrete,
specific theories and even models, making the standard for a ``revolutionary''
change very small indeed, in the limit reaching any kind of conceptual change
whatsoever.


We suggest that there is actually some validity to both moves, that there is a
sort of (weak) self-similarity involved in scientific change.  Every scale of
size and complexity, from local problem solving to big-picture science,
features progress of the ``normal science'' type, punctuated by occasional
revolutions.  For example, in working on an applied research or consulting
problem, one typically will start in a certain direction, then suddenly realize
one was thinking about it wrong, then move forward, and so forth. In a
consulting setting, this reevaluation can happen several times in a couple of
hours.  At a slightly longer time scale, we commonly reassess any approach to
an applied problem after a few months, realizing there was some key feature of
the problem we were misunderstanding, and so forth.  There is a link between
the size and the typical time scales of these changes, with small revolutions
occurring fairly frequently (every few minutes for an exam-type problem), up to
every few decades for a major scientific consensus.  (This is related to but
somewhat different from the recursive subject-matter divisions discussed by
\citealt{Abbott-chaos-of-disciplines}.)  The big changes are more exciting,
even glamorous, but they rest on the hard work of extending the implications of
theories far enough that they can be decisively refuted.

To sum up, our views are much closer to Popper's than to Kuhn's.  The latter
encouraged a close attention to the history of science and to explaining the
process of scientific change, as well as putting on the agenda many genuinely
deep questions, such as when and how scientific fields achieve consensus.
There are even analogies between Kuhn's ideas and what happens in good
data-analytic practice.  Fundamentally, however, we feel that deductive model
checking is central to statistical and scientific progress, and that it is the
threat of such checks that motivates us to perform inferences within complex
models that we know ahead of time to be false.

\section{Why does this matter?}

Philosophy matters to practitioners because they use philosophy to guide their
practice; even those who believe themselves quite exempt from any philosophical
influences are usually the slaves of some defunct methodologist.  The idea of
Bayesian inference as inductive, culminating in the computation of the
posterior probability of scientific hypotheses, has had malign effects on
statistical practice. At best, the inductivist view has encouraged researchers
to fit and compare models without checking them; at worst, theorists have
actively discouraged practitioners from performing model checking because it
does not fit into their framework.

In our hypothetico-deductive view of data analysis, we build a statistical
model out of available parts and drive it as far as it can take us, and then a
little farther.  When the model breaks down, we dissect it and figure out what
went wrong.  For Bayesian models, the most useful way of figuring out how the
model breaks down is through posterior predictive checks, creating simulations
of the data and comparing them to the actual data. The comparison can often be
done visually; see \citet[ch.\ 6]{Gelman-et-al-Bayesian-data-analysis} for
a range of examples.  Once we have an idea about where the problem lies, we can
tinker with the model, or perhaps try a radically new design.  Either way, we
are using deductive reasoning as a tool to get the most out of a model, and we
test the model---it is falsifiable, and when it is consequentially falsified,
we alter or abandon it.  None of this is especially subjective, or at least no
more so than any other kind of scientific inquiry, which likewise requires
choices as to the problem to study, the data to use, the models to employ,
etc.---but these choices are by no means arbitrary whims, uncontrolled by
objective conditions.

Conversely, a problem with the inductive philosophy of Bayesian statistics---in
which science ``learns'' by updating the probabilities that various competing
models are true---is that it assumes that the true model (or, at least, the models among which we will choose or average over) is one of the
possibilities being considered. This does not fit our own experiences of
learning by finding that a model doesn't fit and needing to expand beyond the
existing class of models to fix the problem.

Our methodological suggestions are to construct large models that are capable
of incorporating diverse sources of data, to use Bayesian inference to
summarize uncertainty about parameters in the models, to use graphical model
checks to understand the limitations of the models, and to move forward via
continuous model expansion rather than model selection or discrete model
averaging.  Again, we do not claim any novelty in these ideas, which we have
and others presented in many publications and which reflect decades of
statistical practice, expressed particuarly forcefully in recent times by
\citet{Box-sampling-and-Bayes} and \citet{Jaynes-book}. These ideas, important
as they are, are hardly groundbreaking advances in statistical methodology.
Rather, the point of this paper is to demonstrate that our commonplace (if not
universally accepted) approach to the practice of Bayesian statistics is
compatible with a hypothetico-deductive framework for the philosophy of
science.

We fear that a philosophy of Bayesian statistics as subjective, inductive
inference can encourage a complacency about picking or averaging over existing
models rather than trying to falsify and go further.\footnote{\citet[p.\
  112]{Ghosh-Ramamoorthi} see a similar attitude as discouraging enquiries into
  consistency: ``the prior and the posterior given by Bayes theorem [sic] are
  imperatives arising out of axioms of rational behavior---and since we are
  already rational why worry about one more'' criterion, namely convergence to
  the truth?}  Likelihood and Bayesian inference are powerful, and with great
power comes great responsibility.  Complex models can and should be checked and
falsified.  This is how we can learn from our mistakes.

\subsection*{Acknowledgments}

We thank the National Security Agency for grant H98230-10-1-0184, the
Department of Energy for grant DE-SC0002099, the Institute of Education
Sciences for grants ED-GRANTS-032309-005 and R305D090006-09A, and the National
Science Foundation for grants ATM-0934516, SES-1023176, and SES-1023189.  We
thank Wolfgang Beirl, Chris Genovese, Clark Glymour, Mark Handcock, Jay Kadane,
Rob Kass, Kevin Kelly, Kristina Klinkner, Deborah Mayo, Martina Morris, Scott
Page, Aris Spanos, Erik van Nimwegen, Larry Wasserman, Chris Wiggins, and two
anonymous reviewers for helpful conversations and suggestions.

\bibliography{locusts}
\bibliographystyle{plainnat}

\end{document}